\documentclass[12pt]{amsart}
\usepackage{amsmath, amssymb, amsthm, latexsym}
\usepackage{slashed}
\input xypic
\usepackage{fullpage,amssymb,epic,eepic,epsfig,amsmath}
\newtheorem{theorem}{Theorem}
\newtheorem{proposition}[theorem]{Proposition}
\newtheorem{observation}[theorem]{Observation}

\newtheorem{lemma}[theorem]{Lemma}

\newtheorem{definition}[theorem]{Definition}

\def\Proof{\medskip\noindent{\bf Proof: }}

\def\Z{\mathbb{Z}}

\def\P{\mathbb{P}}

\def\C{\mathbb{C}}
\def\E{\hat{E}}

\def\R{\mathbb{R}}
\def\C{\mathbb{C}}

\def\N{\mathbb{N}}
\def\S{\tilde{S}}
\def\F{\mathbb{F}}

\def\Pi{\mathbb{P}^{\infty}}

\def\qed{\hfill$\square$\medskip}

\def\Zpk{\mathbb{Z}/p^{k}}
\def\Zpk1{\mathbb{Z}/p^{k-1}}

\newcommand{\rref}[1]{(\ref{#1})}

\newcommand{\cform}[3]{\begin{array}{c}
{\scriptstyle #3}\\
#1\\
{\scriptstyle #2}\end{array}}

\newcommand{\fracd}[2]{\frac{\displaystyle #1}{\displaystyle #2}}

\newcommand{\beg}[2]{\begin{equation}\label{#1}#2\end{equation}}
\def\r{\rightarrow}

\def\mh{\mathcal{H}}

\def\F{\mathbb{F}}

\def\S{\mathbb{S}}

\def\sl2{\widetilde{SL_{2}(\Z)}}

\def\sign{\mathrm{sign}}

\def\BZ{\mathbb Z}

\def\a{\alpha}

\def\E{\mathcal E}
\def\P{\mathcal P}
\def\E{\mathcal E}

\def\Pfaf{\operatorname{Pfaf}}

\title[Fermions and Dimers]{On discrete field theory properties of
the dimer and Ising models and their conformal field theory limits}

\author{Igor Kriz}
\address{Igor Kriz\\
Department of Mathematics, \\ University of Michigan
2074 East Hall, 530 Church Street\\
Ann Arbor, MI 48109-1043\\
U.S.A.}
\email{ikriz@umich.edu}
\author{Martin Loebl}
\address{Martin Loebl\\Dept.~of Applied Mathematics and\\
Institute of Theoretical Computer Science (ITI)\\
Charles University \\
Malostranske n. 25 \\
118 00 Praha 1 \\
Czech Republic.}
\email{loebl@kam.mff.cuni.cz}
\author{Petr Somberg}
\address{Petr Somberg\\Mathematical Institute\\
Charles University \\
Sokolovska 83 \\
180 00 Praha 8 \\
Czech Republic.}
\email{somberg@karlin.mff.cuni.cz}

\thanks{
I. Kriz was supported by NSF grant DMS 1102614,
P. Somberg would like to acknowledge a support of the research project MSM 0021620839 and  
the grant GA \v CR 201/08/0397.
}

\begin{document}

\begin{abstract}
 We study various mathematical aspects of discrete models on 
graphs, specifically the Dimer and the Ising models. We focus on 
proving gluing formulas for individual summands of the partition
function. We also obtain partial results regarding conjectured limits realized by 
fermions in rational conformal field theories.
\end{abstract}


\maketitle

\begin{center}


\end{center}

\section{Introduction}

\label{sintro}

The idea that the dimer and Ising statistical models on finite graphs have,
as a limit, the free fermion conformal field theory in dimension $2$
and $1$, respectively, has become a well established theme
in mathematical physics (cf. \cite{onsager}, \cite{dijk}). 
Yet, great challenges remain on the road to developing this into
a mathematically rigorous theory. A part of the
difficulty is in reconciling the languages and concepts
of graph theory and conformal field theory: the structure 
of conformal field theory predicts
many features which are not readily visible on graphs. One must develop mechanisms
how such features originate on the graph level and will fully manifest
themselves in the limit. The purpose of this paper is, in some
sense, to begin this investigation
in earnest. We record some of the phenomena which need to modelled.
We also obtain concrete results as first steps in the desired direction.

\vspace{3mm}
While conformal field theory
was long considered a part of the domain of physics, rigorous mathematical
approaches now do exist. Perhaps most appealing is the ``naive'' approach
due to Graeme Segal \cite{scft}, which asks for a Hilbert space $H$, 
and for each Riemann surface $\Sigma$ with analytically parametrized boundary components,
an element $U(\Sigma)$ (defined up to scalar
multiple) of the Hilbert tensor product of copies of $H$ and its dual indexed
by the boundary components (depending on orientation), satisfying appropriate
``gluing axioms''. The paper \cite{scft} is just an outline, and details
have been since filled in, in part by Kriz and his co-authors \cite{spin,hk,fhk}.

\vspace{3mm}
A competing approach, starting with the notion of vertex algebra, was started by
Borcherds \cite{bor} and Frenkel-Lepowsky-Meurman \cite{flm} (see also
Beilinson-Drinfeld \cite{bd}). This approach
models only a part of the gluing structure postulated by Segal,
but has the advantage that a vertex algebra is a completely algebraic concept, 
which does not require analysis. A vertex algebra, however, models
only a part of conformal field theory, the ``chiral $0$-sector''. 
An extension of the approach which models a whole conformal field theory,
was developed by Huang and Lepowsky \cite{hl}, but it does involve analysis,
although not Hilbert spaces. Efforts to model the entire conformal field theory
as a completely algebraic concept are underway \cite{hkp, kx}.

\vspace{3mm}
The mathematically richest case of conformal field theory is a {\em rational
conformal field theory}, when there is a finite set of {\em labels} $\Lambda$
with an involution $(?)^*$
and the Hilbert space breaks up as
\beg{echiral}{H=\bigoplus_{\lambda\in\Lambda} 
H_{\lambda}\hat{\otimes} \overline{H}_{\lambda^{*}}.}
The Hilbert spaces $H_\lambda$ are called {\em chiral sectors} and are subject
to additional axioms. Parallel to the formula \rref{echiral}, there 
is a formula on graded dimensions, the {\em partition function formula}
(formula \rref{ecftfactor} below).
A major advantage of the vertex algebra approach is that it predicts from 
a vertex algebra alone (i.e. algebraic data) when we get a rational
conformal field theory, and gives an algebraic formula for its sectors
(see also \cite{bd}).

\vspace{3mm}
The free fermion is now completely described mathematically
as a conformal field
theory both in the Segal and vertex algebra formalisms (see \cite{scft, spin,
hl}). Yet, the theory is more subtle than one might think. We describe some
of the issues involved in Section \ref{fcft} below. For one thing, a rational
conformal field theory in the precise sense outlined above actually describes
the {\em bosonisation} of the free fermion. There is a variant of the above
description, called a {\em fermionic rational conformal field theory}, which
is more closely related to the topic of this paper. It is actually the fermionic
rational conformal field theory analogue of the partition function formula
(formula \rref{ecftfactori}, and its analogues for higher
genus) we are interested in. In the case of the free fermion,
this formula \rref{ed3} is precisely a sum over spin structures of a Riemann surface 
of genus $g$, with a factor $1/2^g$, which matches the formula
on graph dimer partition function found by Cimasoni and Reshetikhin \cite{cr1}.

\vspace{3mm}
The present paper is mainly about graphs. How can we model the concepts
mentioned discretely? Fortunately, great strides have already been made. 
Mercat and Kenyon \cite{kast, kenyon} defined critically embedded
graphs on Riemann surfaces, which give discrete models the conformal structure. 
A basic quantum field theory-like
gluing property of the dimer model was observed by Cimasoni and
Reshetikhyn \cite{cr}. 

\vspace{3mm}
A large part of the present paper (Sections \ref{s.dg}, \ref{sec.str}) consists of
examining the gluing property in more detail, both for the Ising and
dimer models. Our main new result in this direction is proving a gluing
formula not just for the entire partition function, but for its
summands, i.e. discrete analogues of \rref{ed3}, \rref{ecftfactori}.
More precisely, we consider {\em arbitrary} labellings, not only critical
ones; this generality, the formula is actually related to quantum field
theory (not conformal field theory). The conformal case, however, is
the one we need to consider for taking scaling limits. For simplicity,
we restrict to the case of a double torus being cut into two surfaces of
genus $1$, but it is straightforward to generalize the formula to more
general cuts along separating curves.

\vspace{3mm}
Our gluing formula is substantially more difficult to obtain than the
formula of \cite{cr}. Instead of cutting along a system of edges of
a graph embedded into a Riemann surface by a
separating curve, we consider a separating curve containing vertices only
of the graph. To obtain a cutting (gluing) formula for partition functions,
we must perform a certain graph-theoretical constrution on the graphs resulting 
from the cut (introducing new subgraphs which we call {\em core} in the
case of the dimer and {\em target} in the case of the Ising model).
Our main gluing theorems are Theorems \ref{thm.d1} and \ref{thm.d2}
of Section \ref{sec.str}.

\vspace{3mm}
To work our way closer toward a limit formula, we need to understand how gluing works 
mathematically in the fermion conformal field theory. This is done in
Section \ref{fermion}. Even though there are both mathematical
and physical approaches (e.g. \cite{w,spin,feingold}), as far
as we know, there is no treatment in the literature which would include
a mathematically precise definition of the data at a closed Riemann
surface in terms of the Dirac operator, and the related cases of
gluing. We write
the Dirac operator on a Riemann surface in a conformally invariant form,
and write down a definition and gluing properties of its determinant
with respect to gluing along curves with antiperiodic spin structure.
(The determinant is the square of the Pfaffian: we restrict to
the case of the determinant because it is much simpler technically.)

\vspace{3mm}
A rigorous
discrete-to-continuous limit formula for one of the chiral partition functions of the
chiral fermion on a torus is the topic of Section \ref{slimit}. Previous work
on this was done by Cohn, Kenyon and Propp who work out a
formula in the case of rectangular
dominoes \cite{ckp} (in fact, in the case of a parallelogram, it was worked
out already by Kasteleyn \cite{kast}, see also Ferdinand
\cite{ferdinand}). A heuristic approach in a more general
situation is done in \cite{dijk}.
A general obstacle to deal with, which we haven't mentioned
above, is the {\em conformal anomaly}, namely the fact that the 
elements $U(\Sigma)$ are defined only up to scalar multiple. 
(In the chiral sector, this propagates into the elaborate structure
of a {\em modular functor}, some of which is described in 
Section \ref{fcft}.)
To deal with this difficulty, one encounters
a phenomenon in nearly all discrete-to-continuous limits in physics:
{\em regularization}. Typically, one encounters a limit which does not
exist, and is ``regularized'' by replacing it with a somewhat different
but compellingly similar expression. In the present setting, we prove for
arbitrary periodic critical embedding in a torus
that the limit does not exist, and can be regularized by replacing,
at one point, a set of parameters by their absolute values, after which
the CFT limit is obtained.

\vspace{3mm}
There is yet another complication in a rigorous approach
to discrete-to-continuous limits 
in the case of the Ising model, which does not arise in the case of
the dimer model: while in the case of the dimer model, the limit
theory is the $2$-dimensional fermion CFT, 
for the Ising model, the limit theory should be the $1$-dimensional fermion
\cite{onsager},
which, despite similar terminology, is considerably
more subtle, and involves a construction of P. Deligne using
super-central simple algebras over $\C$. Whether natural
discrete versions of these concepts exist is an open 
problem, which we discuss at the end of Section \ref{fcft}. 

\vspace{3mm}

\section{Discrete gluing: the set-up}
\label{s.dg}
Let $G= (V(G),E(G))$ be a finite unoriented  graph consisting of vertices $V(G)$
and edges $E(G)$, in which loop-edges and multiple
edges are allowed. We say that the subset $E'\subset E(G)$ is
{\em even} if the graph $(V(G),E')$ has even degree (icluding degree
zero) at each vertex. We say that $M\subset E(G)$ is 
a {\em perfect matching} or {\em dimer arrangement} if the
graph $(V(G),M)$ has degree one at each vertex.  Let $\mathcal{E}(G)$
denote the set of all even subgraphs of $G$, and let $\mathcal {P}(G)$
denote the set of all perfect matchings (dimers) of $G$. 

We assume that an indeterminate $x_e$ is associated with each edge
$e\in E(G)$, and $x_e$ may be evaluated in complex numbers.
 We define the generating functions for even subsets of $E(G)$ resp. for
perfect matchings,  $\E_G$ resp. $\P_G$, as elements of $
\BZ[(x_e)_{e\in E(G)}]$: 
$$\E_G(x)= \sum_{E'\in\mathcal{E}(G) }\,\,\prod_{e\in E'}x_e~,$$
$$\P_G(x)= \sum_{M\in\mathcal{P}(G) }\,\,\prod_{e\in M}x_e~.$$
Polynomial $\P_G(x)$ is also known as the dimer partition function of
graph $G$.
 Knowing the polynomial $\E_G$ is equivalent to knowing the partition
 function 
 $Z^{\mathrm{Ising}}_G$
of the Ising model on
the graph $G$, defined by 
$$
Z_G^{\mathrm{Ising}}(\beta) = Z^{\mathrm{Ising}}_G(x)\Big|
  _{\textstyle{x_e:=e^{\beta J_e} \ \forall  e\in
  E(G)}}
$$
 where the $J_e$ $( e\in E(G))$ are weights (coupling constants)
 associated with the edges of the graph $G$, the parameter $\beta$ is the inverse
 temperature,  and
$$ 
Z_G^{\mathrm{Ising}}(x)= \sum_{\sigma:V(G)\rightarrow \{1,-1\}} \  \prod_{e= \{u,v\}\in E(G)}x_e^{\sigma(u)\sigma(v)}.
$$
The theorem of van der Waerden \cite{vdW} 
states that $Z_G^{\mathrm{Ising}}(x)$ is equivalent to $\E_G(x)$ up to a
change of variables and 
multiplication by a constant factor: 

$$
Z_G^{\mathrm{Ising}}(x)= 2^{|V(G)|}\left(\prod_{e\in E(G)}\frac{x_e+ x_e^{-1}}{2}\right) \E_G(z)\Big|
  _{\textstyle{z_e:= \frac{x_e- x_e^{-1}}{x_e+ x_e^{-1}}}}
$$

Let us assume the vertices of $G$ are numbered
from $1$ to $n$. If $D$ is an orientation of $G$, we denote by $A(G,D)$ the {\em skew-symmetric adjacency
matrix} of $D$ defined as follows: The diagonal entries of  $A(G,D)$ are zero,
and the off-diagonal entries are 
$$A(G,D)_{ij}=\sum \pm x_e~,$$ where the sum is over all edges $e$
connecting vertices $i$ and $j$, and the sign in front of $x_e$ is $1$ if $e$ is
oriented from $i$ to $j$ in the orientation $D$, and $-1$
otherwise. As is well-known, the Pfaffian
of this matrix counts perfect matchings of $G$ {\em with signs:} 
$$\Pfaf A(G,D) = \sum_{M\in \mathcal{P}(G)}\sign(M,D) \prod_{e\in
  M}x_e~,$$ where $\sign(M,D)=\pm 1$. We use this as the definition of the sign of a perfect
matching $M$ with respect to an orientation $D$. 
It is well known that we have 
for any pair of perfect matchings $M,N$
$$
\sign(M,D)\cdot\sign(N,D)= 
(-1)^{\text{\mbox{ the number
of }D-\text{even cycles in} } M\cup N}.
$$ 
An even-length cycle is $D-$even if it has an even number of directed edges in agreement with one way of traversal.

We denote the polynomial 
 $\Pfaf A(G,D)\in \BZ[(x_e)_{e\in E(G)}]$ 
by
$F_D(x)$ and call it the {\em
Pfaffian} associated to the orientation $D$. The following result is
well-known:

\begin{theorem}[Kasteleyn \cite{kast}, Galluccio-Loebl \cite{gl},
  Tesler \cite{tes}, Cimasoni-Reshetikhin \cite{cr1}] \label{Pfaff-formula} If $G$
embeds into an orientable surface of genus $g$, then there exist $4^g$
orientations $D_i$ ($i=1,\ldots, 4^g$) of $G$ such that the perfect matching
polynomial $\P_G(x)$ can be
expressed as a linear combination  of the Pfaffian polynomials
$F_{D_i}(x)$.
\end{theorem}

We call this expression the {\em Arf-invariant formula}, as it is
based on a property of Arf invariant for quadratic forms in
characteristic two. As far as we know, the relationship with Arf
invariant was first observed in \cite{cr1}. 

Let the graph $G$ be embedded in a closed Riemann surface $X$ of genus $g$.
We denote $H:=H_1(X,\F_2)$ the first homology group of $X$ with coefficients in the field $\F_2$.
Recall that $H$ carries a non-degenerate skew-symmetric bilinear form
called $\pmod 2$ {\em intersection form} and denoted by '$\cdot$'.
{Quadratic form on $(H,\cdot)$} associated to  '$\cdot$'
is a function 
$q:H\rightarrow \F_2$ fulfilling
$$
q(x+y)=q(x)+q(y)+x\cdot y
$$ 
for all $x,y \in H$.

 We denote $Q$ the set of quadratic forms on $H$ over $\F_2$, its cardinality being $4^g$. 

Each quadratic form $q\in Q$ determines the {\em signed generating function of even subsets} of $G$
by the formula
\begin{eqnarray}
\label{egen}
\E_{G,q}(x)= \sum_{E' \in \E(G)}(-1)^{q([E'])}\prod_{e\in E'}x_e.
\end{eqnarray}

The generating function $\E_G(x)$ of the even subsets of edges admits analogous {\em Arf invariant formula}, where
$\E_G(x)$ is expressed as a linear combination of its signed modifications 
$\E_{G,q}(x), q\in Q$, see \cite{ls}.

Analogously, partition function of the free fermion conformal field
theory in (spacetime) dimension $2$ on a closed 
Riemann surface $X$ of genus $g$ can be expressed as a linear combination of 
$2^{2g}$ functional Pfaffians of the Dirac operators, each of them corresponding 
to a Spin-structure on $X$. From a physical perspective, this
was noted  in e.g. \cite{agmv}, and rigorous mathematical approaches now also
exist (see Sections \ref{fermion}, \ref{fcft} below).  

Quite recently an extensive effort was invested
in the understanding of the Dimer and Ising models, their criticality resp. asymptotic behaviour
in terms of the conformal invariance and quantum field theory
(conformal field theory in the critical case) of the free fermion and related theories
(see e.g. \cite{agmv}, \cite{difr}, \cite{w}, \cite{cr}, \cite{ls}, \cite{dijk}). 
We continue in this effort and show that the dimer
and the Ising partition functions share with the partition function of the free fermion
a basic {\em gluing} property of their signed modifications. 

We present the gluing formulas for finite graphs embedded in double torus $T^2=T\# T$,
regarded as a connected sum of two $2$-tori (elliptic curves) $T$. 
The restriction is done for the sake of simplicity, i.e. analogous formulas can be 
proven for closed Riemann surfaces of any genus.

\subsection{Exterior algebra of gluing variables}
\label{sub.gluvar}
Let us consider a graph $G$ embedded in a double torus $T^2$.
We denote by $C$ a simple closed curve (a cycle) on $T^2$ with the following properties satisfied:
\begin{enumerate}
\item
 $T^2\setminus C$ has two connected components,  each of which is a torus with a disc bounded by $C$ cut out. 
\item
$C$ intersects embedded graph $G$ in an even number of vertices and in no inner point of an edge.
\item
$G$ has an even number of vertices in each component of $T^2\setminus C$.
\end{enumerate}

\begin{figure}[h]
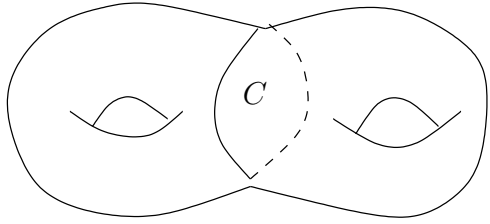

\begin{center}
\input dtor.pstex_t
\end{center}
\caption{The double torus}
\label{fig.dtor}
\end{figure}

Let $v^1,\ldots, v^k$, $k$ even, be the vertices of $G$ on $C$ ordered along $C$. 
We consider the {\em exterior algebra} of $2k$ gluing variables 
$$e_{1}^{1}<e_{2}^{1}<...<e_{1}^{k}<e_{2}^{k}$$
ordered according to the orientation of $C$ and
the choice of a vertex on $C$. Then the vector space on the basis $B$
consisting of all monomials 
\beg{elbas}{e_{\epsilon_1}^{i_1}e_{\epsilon_2}^{i_2}...e_{\epsilon_\ell}^{i_\ell}}
with
$$e_{\epsilon_1}^{i_1}<e_{\epsilon_2}^{i_2}<...<e_{\epsilon_\ell}^{i_\ell}$$
is naturally identified with the exterior algebra
$$\Lambda(e_{1}^{1},e_{2}^{1},...,e_{1}^{k},e_{2}^{k}).$$
The dual exterior algebra  
$$\Lambda(e_{1}^{1},e_{2}^{1},...,e_{1}^{k},e_{2}^{k})^*=
\Lambda(e_{2}^{k*},e_{1}^{k*}, ... , e_{2}^{1*}, e_{1}^{1*})$$
is naturally identified with the vector space on the basis $B^*$ consisting of all monomials 
\beg{elbas1}{e_{\epsilon_\ell}^{i_\ell *}...e_{\epsilon_1}^{i_1 *}.}

Let us consider $Y\in  \Lambda(e_{1}^{1},e_{2}^{1},...,e_{1}^{k},e_{2}^{k})$ written as a linear 
combination of basis elements $Y= \sum_{b\in B}y_bb$, and similarly
$Z\in \Lambda(e_{1}^{1},e_{2}^{1},...,e_{1}^{k},e_{2}^{k})^*$ written as $Z= \sum_{b'\in B^*}z_{b'}b'$.
Then the scalar product  $<Y,Z>$ is defined by
$$
<Y,Z>= \sum_{(b,b')}<y_bb,z_{b'}b'>,
$$
where
$<y_bb,z_{b'}b'>$ is written in the following way. Let $b= e_{\epsilon_1}^{i_1}e_{\epsilon_2}^{i_2}...e_{\epsilon_\ell}^{i_\ell}$
and $b'= e_{\a_m}^{j_m*}...e_{\a_1}^{j_1*}$. Then $<b,{b'}>= 0$ unless
$m= \ell$ and $\epsilon_p= \a_p, i_p= j_p$ for each $p$, in which case $<y_bb,z_{b'}b'>= y_bz_{b'}$.

We first illustrate the gluing by writing down a gluing formula for the unsigned partition function $\E_G(x)$. 
This was done in a different way for the dimer partition function $\P_G(x)$  in \cite{cr}.

\subsection{Discrete unsigned gluing} 
\label{sub.uns}
We denote by $G_1, G_2$ the subgraphs of $G$ supported on connected components of $T^2\setminus C$ along with 
the vertices of $G$ supported on $C$. Hence, if $v$ is a vertex of $G$
embedded in $C$ then the set $B_v$ of the edges incident with $v$ is partitioned into two sets $B_1(v), B_2(v)$,
the first belonging to $G_1$ and the second to $G_2$.

\begin{figure}[h]
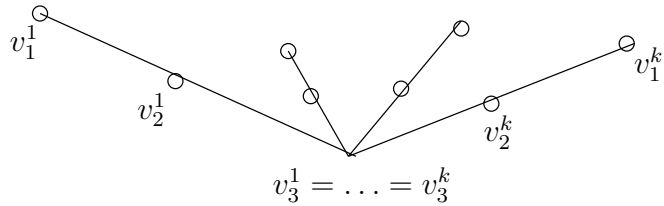

\begin{center}
\input center.pstex_t
\end{center}
\caption{The target $S$}
\label{fig.center}
\end{figure}

We construct an auxiliary graph, called the {\em target}, as follows.
For each $i=1,\ldots, k$, let us introduce the path $L(i)$ of length 2 containing the vertices 
 $v^i_1, v^i_2,v^i_3$ and the edges $\{v^i_3 v^i_2\}, \{v^i_2 v^i_1\}$.  The target $S$ is obtained 
 as the union of paths 
 $L(i), i=1,\ldots, k$ by identifying all $k$ vertices $v_3^i, i=1, \ldots, k$ into one single vertex 
 (see Figure \ref{fig.center} for the target $S$). The target $S$ allows to construct the graph $G^S_j$ from 
 $G_j$ and $S$ by identifying, for each $i= 1, \ldots, k$ the vertices $v^i$ and $v_1^i$.  

 The edge-variables of $S$ are the {\em gluing variables}. We denote by
$y_{1v^i_3v^i_2}= e^i_1$, $y_{1v^i_2v^i_1}= e^i_2$ the gluing variables in $G^S_1$ and by
$y_{2v^i_3v^i_2}= e^{i*}_2$, $y_{2v^i_2v^i_1}= e^{i*}_1$ the gluing variables in $G^S_2$. 

The vector $y_1$ of variables of $G^S_1$ consists of $(x_e)_{e\in E(G)}$ and the gluing variables
 $(e_{1}^{1},e_{2}^{1},...,e_{1}^{k},e_{2}^{k})$. The vector $y_2$ of variables of $G^S_2$  
 consists of $(x_e)_{e\in E(G)}$ and the gluing variables 
$(e_{2}^{k*},e_{1}^{k*}, ... , e_{2}^{1*}, e_{1}^{1*})$.
The following theorem simply follows from the definition of scalar product in the exterior algebra.

\begin{theorem}
\label{thm.uns}
$$
\E_G(x)= <\E_{G^S_1}(y_1), \E_{G^S_2}(y_2)>.
$$
\end{theorem}

\section{Discrete signed gluing}
\label{sec.str}
   
Many results in this section exploit the structure of $g$-graph,
extensively used and explained in detail in, e.g. \cite{LM}, \cite{ls} 
and references therein.
\begin{definition}
\label{def.highway}
The {\em highway surface} $\S_g$ consists of a base polygon $R_0$ and 
collection of bridges $R_1, \ldots, R_{2g}$, where 
\begin{enumerate}
 \item 
$R_0$ is a convex $4g-$gon with vertices $a_1, \ldots, a_{4g}$ numbered clockwise,
\item
Each $R_{2i-1}$ is a 
rectangle
with vertices $x(i,1), \ldots, x(i,4)$ numbered clockwise, and glued
to $R_0$ in a way that its edge $[x(i,1),x(i,2)]$ is identified with the edge
 $[a_{4(i-1)+1},a_{4(i-1)+2}]$
and the edge $[x(i,3),x(i,4)]$ is identified with the edge $[a_{4(i-1)+3},a_{4(i-1)+4}]$,
\item
Each $R_{2i}$ is a 
rectangle
with vertices $y(i,1), \ldots, y(i,4)$ numbered clockwise. It is glued
to $R_0$ so that its edge $[y(i,1),y(i,2)]$ is identified with the edge $[a_{4(i-1)+2},a_{4(i-1)+3}]$
and the edge $[y(i,3),y(i,4)]$ is identified with the edge
$[a_{4(i-1)+4},a_{4(i-1)+5}]$. Notice that indices are considered modulo $4g$.
\end{enumerate}
\end{definition}

There is an orientation-preserving map $\Phi$ of $\S_g$
into the plane $\R^2$, which is an immersion up to the images of the bridges
$R_{2i}$ and $R_{2i-1}$ intersecting in a square for each $i=1, \ldots
g$. 

Let us assume the graph $G$ is embedded into a closed orientable surface $X$
of genus $g$. The combinatorial model of $X$ is $\S_g$ unified with an additional disk
$R_\infty$ glued to the boundary of $\S_g$. By an isotopy property of the embedding,
we may assume that $G$ does not meet the disk $R_\infty$ and 
moreover, all
vertices of $G$ lie in the interior of $R_0$. We may also assume that
the intersection of $G$ with any of the rectangular bridges $R_i$
consists of disjoint straight
lines connecting the two sides 
of $R_i$ glued to the base polygon $R_0$.
This assumption is not really necessary, but makes all considerations below more transparent.  
The composition of the embedding of $G$
into $\S_g$
with the immersion $\Phi$ yields a drawing $\varphi$ of
$G$ in the plane $\R^2$. A planar drawing of $G$ obtained in this way will be
called {\em special}.
Observe that double points of a
special drawing can only come from the intersection of the images of
bridges under the immersion $\Phi$ of $\S_g$ into the plane. Thus every double
point of a special drawing lies in one of the squares
$\Phi(R_{2i})\cap \Phi(R_{2i-1})$.

In what follows we employ this machinary to produce gluing formulas.

\subsection{The Ising model}
\label{sub.is}
As for the notation and basic results, we follow \cite{ls}.
\begin{definition}
\label{def.rrr}
Let $G$ be a graph embedded in $\S_g$ and let $e$ be an edge of $G$. By definition, the embedding of $e$ intersects each bridge $R_i$  
in a collection of disjoint straight lines, whose number is denoted $r_i(e)$.   For a subset $A$ of edges of $G$ we denote $r(A)$ 
the vector of length $2g$, whose $i$-th component is
$r(A)_i= \sum_{e\in A}r_i(e)$.
\end{definition}

We observe that two even subsets of edges $A,B$ belong to the same homology class 
in $H$ if and only if $r(A)= r(B)\pmod 2$.  We pick a basis 
$a_1, b_1, \ldots , a_g,b_g$ of $H$, where $a_i, b_i$ correspond to the even subsets $A_i, B_i$ satisfying $r(A_i)_{2i-1}= r(B_i)_{2i}= 1$ and all remaining components of $r(A_i)$ and $r(B_i)$  are zero. Note that each quadratic form in $Q$ is uniquelly determined by its values on any basis 
of the underlying vector space. 

The following Theorem is a basic result of (\cite{ls}). 
\begin{theorem}
\label{thm.is}
Let $G$ be a graph embedded in a closed Riemann surface $X$ of genus $g$ and let $q$ be a quadratic form on the homology group $H$. 
Let $E'$ be an even subset of edges of $G$ and let $r= r(E')$. Then 
$$
(-1)^{q([E'])}= (-1)^{m(E')},
$$
where
$$
m(E')= \sum_{i=1}^{g}r_{2i-1}r_{2i}+ r_{2i-1}q(a_i)+ r_{2i}q(b_i).
$$
\end{theorem}

The importance of Theorem \ref{thm.is} is related to the following observation. 
Let us consider an embedding of $G$ in the highway surface $\S_2$, for simplicity considered of genus $g=2$. 
A quadratic form $q\in Q$ naturally determines two quadratic forms $q_1, q_2$ on the torus by 
$q_1(a_1)= q(a_1), q_1(b_1)= q(b_1)$ resp. $q_2(a_1)= q(a_2), q_2(b_1)= q(b_2)$ together with 
trivial values on the other basis elements not written explicitly. 

The signed version of Theorem \ref{thm.uns} is of the following form.

\begin{theorem}
\label{thm.d2}
We have
\begin{eqnarray}
\E(G,q,x)= <\E(G^S_1,q_1,y_1), \E(G^S_2,q_2,y_2)>.
\end{eqnarray}
\end{theorem}
\Proof
Using Theorem \ref{thm.uns} the proof is easily reduced to 

{\bf Claim.}
Let $E'$ be an even subset of edges of $G$. For $i=1,2$, there is exactly one extension  $E'_i$ of $E'\cap G_i$ to an even subset 
of $G^S_i$ by edges in the target. 
Moreover, for arbitrary $q\in Q$, 
$$
(-1)^{q([E'])}= (-1)^{q_1([E'_1])}(-1)^{q_2([E'_2])}.
$$
The Claim  follows directly from Theorem \ref{thm.is} and definitions of $q_1, q_2$. This finishes the proof
of Theorem \ref{thm.d2}.

\qed

\subsection{The dimer model}
\label{sub.dim}
In this subsection we present a gluing formula for the signed dimer partition function of a graph $G$ embedded in the closed Riemann surface of genus $2$.
Analogously to the case of the Ising partition function, we first introduce an auxiliary graph called the {\em core} and denoted $S'$. The construction of the core follows the procedure of reduction of the Ising model to the dimer model, see \cite{f}.

We recall that in Subsection \ref{sub.uns}
we denoted by $G_1, G_2$ the subgraphs of $G$ supported on connected components of $T^2\setminus C$ together with 
the vertices of $G$ supported on $C$. If $v$ is a vertex of $G$
embedded in $C$ then the set $B_v$ of the edges incident with $v$ is partitioned into two sets $B_1(v), B_2(v)$,
the first belonging to $G_1$ and the second to $G_2$.
For each $i=1,\ldots, k$, $L(i)$ denotes the path of length 2 consisting of the vertices 
 $v^i_1, v^i_2,v^i_3$ and the edges $\{v^i_3 v^i_2\}, \{v^i_2 v^i_1\}$.  

\begin{figure}[h]
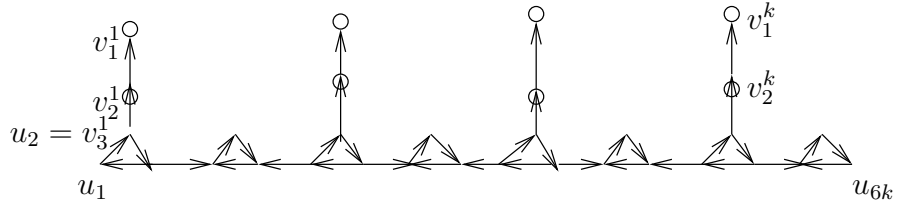

\begin{center}
\input core.pstex_t
\end{center}
\caption{The orientation of the core $S'$}
\label{fig.core}
\end{figure}

To construct the core $S'$, we first consider the path consisting of $6k$ new vertices $u_1,\ldots, u_{6k}$ and edges 
$\{u_i, u_{i+1}\}$, $i=1, \ldots 6k-1$.
Then we add the edges $\{u_{3j-2},u_{3j}\}, j = 1,...,2k$ and identify, for $i=1,\ldots, k$, the
 vertex $u_{6i-4}$ with $v_3^i$ (see Figure \ref{fig.core} for an orientation of the core).

Finally, for $j=1,2$ we construct the graphs $G^{S'}_j$ from the union of $G_j$ and $S'$ by identifying, for each $i= 1, \ldots, k$ 
the vertices $v^i$ and $v_1^i$.  

The relevance of the core stems from the following observation, see e.g. \cite{f}.

\begin{observation}
\label{o.deltt}
Let $A\subset \{v_1^i, i=1, \ldots, k \}$. Then $S'$ has a matching completely covering all vertices of $S'\setminus A$ if and only if $|A|$ is even. Moreover, if  $|A|$ is even, the matching completely covering
 $S'\setminus A$ is uniquelly determined.
\end{observation}

 The edge-variables of the paths $L(i)$ are, as in the construction of the target, the {\em gluing variables}. For $i= 1,\ldots, k$, we denote by 
$z_{1v^i_3v^i_2}= e^i_1$, $z_{1v^i_2v^i_1}= e^i_2$ the gluing variables in $G^{S'}_1$ resp.
$z_{2v^i_3v^i_2}= e^{i*}_2$, $z_{2v^i_2v^i_1}= e^{i*}_1$ the gluing variables in $G^{S'}_2$.

We let vector $z_1$ of variables of $G^{S'}_1$ consist of $(x_e)_{e\in E(G)}$, the gluing variables
 $(e_{1}^{1},e_{2}^{1},...,e_{1}^{k},e_{2}^{k})$ and the remaining edge variables be set to $1$. We let vector $z_2$ of the variables of $G^{S'}_2$  
 consist of $(x_e)_{e\in E(G)}$, the gluing variables $(e_{2}^{k*},e_{1}^{k*}, ... , e_{2}^{1*}, e_{1}^{1*})$, and 
 the remaining edge variables be set to $1$.

We are now ready to discuss gluing formula for unsigned partition function.
The following Theorem simply follows from the definition of the trace in the exterior algebra, using Observation \ref{o.deltt}.

\begin{theorem}
\label{thm.uns1}
$$
\P_G(x)= <\P_{G^{S'}_1}(z_1), \P_{G^{S'}_2}(z_2)>
$$
\end{theorem}

In the next step we consider the signed case. We first modify $G$ by adding the edges $\{v^i, v^{i+1}\}, i=1, \ldots, k-1$ 
and let the variables associated to these
edges be equal to $0$. Hence, adding this set of edges does not change the signed dimer partition function of $G$.
Let us denote the modified graph by $G'$.

A result of \cite{tes} (see also recent \cite{LM}) can be summarized, in the case of the double torus, as follows: 

\begin{theorem}
\label{thm.pfaf}
Let $G$ be a graph embedded in a highway surface $\S_2$.
Let $D$ be an orientation of $G$ so that for each bridge $R$, the immersion of $R_0\cup R$
into the plane has the property that each inner face has an odd number of edges oriented
 clokwise in $D$. Let $M$ be a perfect matching of $G$. Then 
$$
\sign(M,D)= (-1)^{c(M)},
$$
 where $c(M)$ denotes
the number of self-intersections among edges of $M$ inside the special drawing of $G$
in  the plane.
\end{theorem} 

Let $G$ be the graph embedded in the double torus $T^2$. We consider the embedding of $G'$ in the highway surface 
$\S_2$ determined by an embedding of $G$ in the closed Riemann surface $T^2$ of genus $2$ . 
Note that the new edges attached among the vertices of $C$ form a chord of the base polygon $R_0$,
which separates the two pairs of bridges. 
We recall that $v^1, \ldots, v^k$ is the order of the vertices of $G$ on $C$ along this chord, and that $k$ is even. 

 It is not difficult to observe that there is an orientation $D_0$ of $G'$
fulfilling the properties of Theorem \ref{thm.pfaf}, where the orientations of the edges $\{v^i, v^{i+1}\}, i=1, \ldots, k-1$, 
form a directed path of odd length from $v^1$ to $v^k$.

Let $D$ be an orientation of $G$.
We denote by $D_1$ the restriction of $D$ to $G_1$ and by $D_2$
the restriction of $D$ to $G_2$. 
 We consider the orientation of the core $S'$, indicated on Figure \ref{fig.core}. 
 By slight abuse of notation we denote it also by $S'$.  If $D$ is an orientation of 
 $G$ then we denote by $D^{S'}_1, D^{S'}_2$ the orientations of
$G^{S'}_1, G^{S'}_2$ induced by $D$ and $S'$.

\begin{theorem}
\label{thm.d1}
Let $D$ be an orientation of $G$. Then
$$
F_D(x)= <F_{D^{S'}_1}(z_1), F_{D^{S'}_2}(z_2)>.
$$
\end{theorem}

\Proof
Let $M$ be a perfect matching of $G$. It follows from Observation \ref{o.deltt} that for $i=1,2$, there is exactly one extension  $M_i$ of $M\cap G_i$ to perfect 
matching of $G^{S'}_i$ by the core edges. Theorem \ref{thm.d1} follows from Theorem \ref{thm.uns1} and 

{\bf Claim.}
$$
\sign(M,D)= \sign(M_1, D^{S'}_1)\cdot\sign(M_2, D^{S'}_2).
$$
 We first show that the Claim  holds for the orientation $D_0$ introduced above. 
Since $M$ is a subset of edges of $G$, we consider $D_0$ restricted to $G$ only.
We have by Theorem \ref{thm.pfaf}
$$
\sign(M,D_0)= (-1)^{c(M)},
$$
 where $c(M)$ denotes
the number of self-intersections among the edges of $M$ in the special drawing of $G'$
in  the plane. From the construction of $G'$ in Subsection \ref{sub.dim} we have
$$
c(M)= c_1(M)+ c_2(M),
$$
where $c_1(M)$ denotes
the number of self-intersections among the edges of $M$ that belong to the first pair of the bridges of the highway surface $\S_2$, and $c_2(M)$ denotes
the number of self-intersections among the edges of $M$ that belong to the second pair of the bridges of the highway surface $\S_2$.

Finally, the orientations $D_{1}^{S'}, D_{2}^{S'}$ are constructed from $D_0$ in such a way that they also
 satisfy the properties required by Theorem \ref{thm.pfaf}. Hence
\begin{eqnarray}
& & \sign(M_1,D^{S'}_{1})= (-1)^{c_1(M_1)},
\nonumber \\
& & \sign(M_2,D^{S'}_{2})= (-1)^{c_2(M_2)}.
\end{eqnarray}
Summarising, the Claim holds for the orientation $D_0$. An arbitrary orientation $D$ of $G$ is obtained from $D_0$
by reversing orientation in a subset $S(D)$ of the directed edges of $D_0$, so
$$
\sign(M,D)= \sign(M,D_0)(-1)^{|S(D)\cap M|} 
$$
and
\begin{eqnarray}
& & \sign(M_1, D^{S'}_1)= \sign(M_1, D_{1}^{S'})(-1)^{|S(D)\cap M_1|}, \nonumber \\ 
& & \sign(M_2, D^{S'}_2)= \sign(M_2, D_{2}^{S'})(-1)^{|S(D)\cap M_2|}.
\end{eqnarray}
This finishes the proof of the last Claim and at the same time of Theorem \ref{thm.d1}.

\qed


\section{Fermions on compact surfaces I}

\label{fermion}

\vspace{3mm}
In this section, we present some
facts about the free fermion conformal field theory in the context
of the desired limit from the discrete to
the continuous case. From the physics point of view, this material is more than
two decades old, see e.g. \cite{agmv,w}. A mathematical approach
was outlined in \cite{scft}, and further developed in \cite{spin, hk}.
Nevertheless, the mathematical literature as it stands
is deficient in an important point (explained below), which happens to be
the key to 
addressing the limit questions
we are interested in. Treating this gap is the main purpose of the present section.

\vspace{3mm}
The point is that in the mathematical treatments of CFT, {\em closed} Riemann surfaces
(or compact Riemann surfaces without boundary)
play a special role, and a technical complication. 
It is a common technique, in fact, to consider {\em only} (compact) Riemann surfaces
without closed connected components, and argue that data at closed surfaces can
be recovered by cutting out a holomorphic disk, and giving compatible isomorphisms of the
data as the disk varies. This is the approach taken in \cite{scft,spin}.

\vspace{3mm}
For our purposes, however, this is unsatisfactory. The reason is that from
the physical point of view, the case of a closed surface is in fact more fundamental
(\cite{w}), because it is where contact with quantum field theory is made.
The data we are supposed to have in this case is the determinant (resp. Pfaffian)
of the Dirac
operator, which is analogous to the discrete data we studied in sections
\ref{s.dg}, \ref{sec.str}.

\vspace{3mm}
The purpose of this section, then, is to explain how to define the determinant
precisely, and to describe gluing from the data on a surface with boundary as described
in \cite{scft,spin} to the Dirac operator data in the closed surface case.
For simplicity, we treat only the determinant, which corresponds to the free fermion
of spacetime dimension $2$, and not the Pfaffian, which is more complicated.
Also for reasons of simplicity, we discuss only antiperiodic boundary conditions.
The role of all spin structures will be reviewed in more detail in Section \ref{fcft}
below. A discrete analogue of the `boundary to no boundary' gluing case will be,
in fact, exhibited in Section \ref{slimit} below, where we give one prototype
of a rigorous regularized CFT limit computation.

\vspace{3mm}
Let $\Sigma$ be a compact Riemann surface. A spin structure on $\Sigma$
can be identified with a choice of a holomorphic line bundle
$\underline{\Omega^{1/2}}$ on $\Sigma$ and an isomorphism
$$\underline{\Omega^{1/2}}\otimes
\underline{\Omega^{1/2}}\cong \underline{\Omega^{1}}.$$
(We identify a holomorphic bundle with its sheaf of sections $\underline{M}$;
we denote by $M$ the space of global holomorphic sections; $\underline{M}_{sm},
M_{sm}$ denotes the corresponding notions in the smooth category.)
The bundles of spinors on $\Sigma$ are $\underline{\Delta^-}=\underline{\Omega^{-1/2}}$
and $\underline{\Delta^-}=\underline{\Omega^{+1/2}}$. 

The Dirac operator
$$\Delta^{-}_{sm}\r \Delta^{+}_{sm}$$
ordinarily depends on metric (i.e. is only defined on a spin Riemann manifold),
which is undesirable for the purposes of CFT. However,
a modified conformally invariant construction can be obtained
as follows: We define a pairing
\beg{ed1}{\slashed{D}:\Omega^{1/2}_{sm}\otimes_\C \Omega^{1/2}_{sm}\r \C
}
in a local holomorphic coordinate $z$ as follows:
\beg{ed2}{fdz^{1/2}\otimes gdz^{1/2}\mapsto\int_\Sigma\frac{\partial f}{\partial 
\overline{z}}gdzd\overline{z}.
}
To see that \rref{ed2} is invariant under holomorphic
coordinate change, let $z=z(y)$. Using the $y$ coordinate,
we compute
$$
fdz^{1/2}\otimes gdz^{1/2} = f\left(\frac{dz}{dy}\right)^{1/2}dy^{1/2}
\otimes g \left(\frac{dz}{dy}\right)^{1/2}dy^{1/2}\mapsto$$ $$
\mapsto \int \frac{\partial f}{\partial 
\overline{y}}g\frac{dz}{dy}dyd\overline{y}=\int\frac{\partial f}{\partial 
\overline{z}}g\frac{dz}{dy}\frac{d\overline{z}}{d\overline{y}}dyd\overline{y}=
\int\frac{\partial f}{\partial 
\overline{z}}gdzd\overline{z}.
$$
Also note that by Stokes' theorem, $\slashed{D}$ is an antisymmetric pairing,
since
$$\left(\frac{\partial f}{\partial\overline{z}}g+f\frac{\partial g}{\partial\overline{z}}
\right)dzd\overline{z}=d(fgdz).$$
One defines (up to scalar multiple) the partition function of the
$1$-dimensional chiral fermion on $\Sigma$ as
$$pf(\slashed{D}).$$
The partition function of a fermion with both chiralities
is
\beg{ed3}{\frac{1}{2^g}\cform{\sum}{\sigma}{}\alpha(\sigma) pf\left(
\begin{array}{cc}
\slashed{D}_\sigma & 0\\
0 & \widetilde{\slashed{D}}_\sigma
\end{array}
\right)
}
where $\widetilde{\slashed{D}}$ is the construction analogous to $\slashed{D}$
with ``holomorphic'' replaced by ``anti-holomorphic'' 
and $\alpha(\sigma)$ is the Arf invariant. 

To define the Pfaffian, one proceeds as follows. For an antisymmetric 
Fredholm operator
$$\Phi:H\r H^*$$
where $H$ is a Hilbert space (we are not in this situation, but may
achieve it by suitable Hilbert-completion), there is a {\em Pfaffian line}
$Pf(\Phi)$ (well defined up to canonical isomorphism), and a Pfaffian element
\beg{ed4}{pf(\Phi)\in Pf(\Phi)
}
(see \cite{ps}, Chapter 12).
For simpicity, we will describe only the square of \rref{ed4}, 
\beg{ed5}{det(\Phi)\in Det(\Phi)\;\;\text{for $\Phi:H_1\r H_2$ Fredholm of index $0$.}
}
Choose an invertible operator $J:H_1\r H_2$ such that
$$T=J-\Phi$$
is trace class. (For any compact operator $K:H_1\r H_2$,
there exist unique numbers $s_1\geq s_2\geq...\geq 0$, $\lim s_n=0$
and bases $(e_n)$, $(f_n)$ of $H_1,H_2$, respectively, such that
$K(e_n)=s_n f_n$. The numbers $s_n$ are called {\em singular values}.
$K$ is {\em trace class} (resp. {\em Hilbert-Schmidt})
if $\sum |s_n|<\infty$ (resp. $\sum (s_n)^2<\infty)$.)

Now one defines
$$Det(\Phi):=\langle (J)\rangle$$
(=the line with free formal generator $(J)$), and
$$det(\Phi):=det(I-TJ^{-1})\cdot(J)\in Det(\Phi)$$
(the determinant is defined as a number because
$TJ^{-1}$ is trace class). For another choice
$$\Phi=J^\prime-T^\prime$$
with $J^\prime$ invertible and $T^\prime$ trace class,
we have
$$J^\prime J^{-1}=I+(T^\prime-T)J^{-1},$$
so the canonical iso
$$\langle(J^\prime)\rangle\r \langle(J)\rangle$$
is
$$(J^\prime)\mapsto (J)\cdot det(I+(T^\prime-T)J^{-1}).$$


\vspace{3mm}

One has $det(\Phi)=0$ unless $\Phi$ is invertible. To define
$det(\Phi)$ as a number, one must specify a generator of
$Det(\Phi)$. This is called {\em regularization}. 

Let us turn back to Pfaffians.
In our situation, this amounts to picking a section of
the pullback of the bundle $Pf(\slashed{D}_\Sigma)$ on the 
Teichm\"{u}ller space. This can be done using the structure
of CFT.

For a Riemann surface $\Sigma$ with spin structure 
with analytically parametrized boundary
such that the induced spin structure on each boundary component
is antiperiodic, there is a non-degenerate quadratic form
on the space $\Omega^{1/2}_{\partial \Sigma}$ of complex-valued
$1/2$-forms on the boundary given by
$$B(\omega,\eta)=\int_{\partial\Sigma}\omega\eta.$$
For each boundary component, $1/2$-forms which extend to the
unit disk holomorphically form a maximal isotropic subspace
$(\Omega^{1/2}_{S^1})^-$ of $(\Omega^{1/2}_{S^1})$. One defines
\beg{ed+}{\mh_{S^1}=(\Lambda(\Omega^{1/2}_{s^1})^-)^\wedge
}
(here $?^\wedge$ denotes Hilbert completion,) and
$$\mh_{\partial \Sigma}:=\widehat{\bigotimes} (\mh_{S^1})^*\hat{\otimes} 
\widehat{\bigotimes} \mh_{S^1}$$
where $\hat{\otimes}$ denotes Hilbert tensor product,
and the tensor products on the right hand side are over
boundary components with inbound resp. outbound orientation (the
orientation of the boundary component of the standard disk $D$ parametrized
by the identity is outbound).

In fact, more explanation is in order: In the present situation, note that
the symmetric bilinear form $B$ together with the Real structure (complex 
conjugation) specifies an inner product on $\Omega^{1/2}_{S^1}$
(hence on $\Omega^{1/2}_{\partial \Sigma}$) which is the specific one
we use in forming the Hilbert completion, (and hence also defining,
say, the inner product on $\Lambda((\Omega^{1/2}_{S^1})^-)$ etc.

Now the inclusion
\beg{edr}{(\Omega^{1/2}_{\Sigma})^\wedge\subset (\Omega^{1/2}_{\partial\Sigma})^\wedge
}
where the second map is orthogonal projection is an inclusion of a maximal
isotropic subspace which is {\em restricted} in the sense that its orthogonal
projection to $((\Omega^{1/2}_{\partial\Sigma})^-)^\wedge$ is Fredholm,
while the orthogonal projection to the orthogonal complement
$((\Omega^{1/2}_{\partial\Sigma})^+)^\wedge$ is Hilbert-Schmidt.

In this situation, $\Omega^{1/2}_{\Sigma}$ specifies a line $Pf(\Omega^{1/2}_{\Sigma})$
and a linear map
\beg{ed+++}{Pf(\Omega^{1/2}_{\Sigma})\r\mh_{\partial \Sigma}.
}
Now when $\Sigma$ is a closed Riemann surface with spin structure, define
\beg{ednew*}{Pf(\Omega^{1/2}_{\Sigma}):=Pf(\slashed{D}_\Sigma).}

\begin{lemma}
\label{lg}
Let $\check{\Sigma}$ be a Riemann surface obtained from
$\Sigma$ by gluing some inbound boundary components to
outbound boundary components in a way so as to preserve
parametrization. Then this data specifies an isomorphism
\beg{edi}{\diagram Pf(\Omega^{1/2}_{\Sigma})\rto^\cong_\iota &
Pf(\Omega^{1/2}_{\check{\Sigma}})\enddiagram
}
and a commutative diagram where the rows are given by CFT 
structure
\beg{ednew+}{\diagram
Pf(\Omega^{1/2}_{\Sigma})\rto\dto_\iota & \mh_{\partial\Sigma}\dto^{tr}\\
Pf(\Omega^{1/2}_{\check{\Sigma}})\rto & \mh_{\partial \check{\Sigma}}.
\enddiagram
}
(In the case when $\check{\Sigma}$ is closed, the bottom right
corner is $\C$, and the bottom row is defined to be
$pf(\slashed{D})^{-1}$.)
\end{lemma}

We will show only the construction of the map \rref{edi} in the case
when $\check{\Sigma}$ is closed (i.e. where the definitions of the two
Pfaffian lines are different). Also, for simplicity, we will 
only treat the square, i.e. the analogous construction for determinant
lines: the construction for Pfaffians is similar but more elaborate.

Denote by $(\Omega^{1/2}_{0})_\Sigma$ the subspace of the space
of smooth sections $\omega$ of $(\underline{\Omega}^{1/2}_{sm})_\Sigma$
with the property that for any two boundary components in $\partial \Sigma$
which are to be glued in $\check{\Sigma}$,
$$\overline{\partial}\omega|_{c_1}=\overline{\partial}\omega|_{c_2}.$$
Denote the image of $\partial \Sigma$ in $\check{\Sigma}$ by
$\partial\Sigma/2$. Then we have a commutative diagram
\beg{edii}{\diagram
&\Omega^{1/2}_{\Sigma}\dto_\subset \rto &\Omega^{1/2}_{\partial \Sigma/2}\dto^\cong\\
(\Omega^{1/2}_{sm})_{\check{\Sigma}}\dto^{\slashed{D}}\rto &
(\Omega^{1/2}_{0})_{\Sigma}\dto^{\slashed{D}}\rto &
\Omega^{1/2}_{\partial \Sigma/2}\\
((\Omega^{1/2}_{sm})_{\check{\Sigma}})^*\rto^= &((\Omega^{1/2}_{sm})_{\check{\Sigma}})^*&
\enddiagram
}
where all rows and columns are short exact with the exception of the top row
and left column, which are Fredholm (to reduce clutter, we omit the symbols for
Hilbert completions). The top row differs from \rref{edr} by a trace class operator
(note: ultimately, this difference gives rise to the eta-factor!), so it suffices
to prove that for any diagram of Hilbert spaces and bounded operators
\beg{ediii}{\diagram
& D\rto^G\dto^j &B\dto^=\\
C\dto_F\rto^i & W\dto\rto & B\\
A\rto^= &A &
\enddiagram
}
where all the rows and columns are short exact except $F$ and $G$ which
are Fredholm, we have a canonical isomorphism
\beg{ednewsq}{Det(F)\cong Det(G).}
But note that in \rref{ediii}, we may write $B=W/C$, $A=W/D$,
so both sides are obviously canonically isomorphic to
\beg{ednewstar}{Det(i\oplus j:C\oplus D\r W).}
The commutativity of the square \rref{ednew+} is proved in \cite{spin}.
\qed

Now for the case $\Sigma =A_q$ (the standard annulus in $\C$ with boundary components
parametrized by $z$, $qz$ and antiperiodic spin structure), the
Lemma specifies a canonical isomorphism
$$Pf(\Omega^{1/2}_{A_q})\otimes Pf(\Omega^{1/2}_{D})\r Pf(\Omega^{1/2}_{D}),$$
which specifies a canonical element
\beg{edv}{\iota_q\in Pf(\Omega^{1/2}_{A_q}).
}
Considering the elliptic curve $\check{\Sigma}=E_q$ over $\C$,
the image
$$\xi_q:=\iota(\iota_q)\in Pf(\Omega^{1/2}_{E_q})$$
is the desired regularization. With respect to this regularization, the partition
function is one of the two Jacobi theta functions involving odd powers of $q^{1/2}$
(depending on the spin of the gluing), times the factor
\beg{edvi}{\cform{\prod}{n>0}{}(1-q^n)=\eta(\tau)q^{1/24},\; q=e^{2\pi i\tau}.
}
The remaining Arf invariant $0$ theta function is not accessible in this way,
but may be obtained by cutting around a different simple curve in $E_q$
of anti-periodic spin structure; after multiplying by the factor $q^{-1/24}$,
it is related to the remaining Arf invariant $0$ partition functions by
a modular transformation. The Arf invariant $1$ partition function is $0$.

\vspace{3mm}

In the next section, we will see that the partition
function of the chiral fermion is, in a special case where the
boundary can be easily modelled discretely, a ``regularized limit''
of partition functions of dimer models on critically embedded graphs, using
a discrete analogue of Lemma \ref{lg}.

In Section \ref{s.dg} Theorem \ref{thm.uns1}, we exhibited a discrete analogue of Lemma \ref{lg}
for arbitrary (not necessarily critically embedded) finite graphs, embedded into a 
Riemann surface of genus $2$.

\vspace{3mm}

\section{A limit formula for the discrete dimer model on a torus}
\label{slimit}

\vspace{3mm}
We use here the setup of Kenyon \cite{kenyon}. Consider a bipartite planar graph
with vertices colored black and white. The set of black (resp. white) vertices is
denoted by $B$ resp. $W$. A critical embedding into a plane is such that the vertices
of each face lie on a circle, and all these circles have the same
radius. One defines a square matrix $\overline{\partial}$
with rows and
columns indexed by $B\cup W$ as follows: If $v_1,v_2$ are not adjacent, then
$\overline{\partial}(v_1,v_2)=0$. If $w$ and $b$ are adjacent vertices, $w$ being
white and $b$ being black, then
\beg{eanti}{\overline{\partial}(w,b)=-\overline{\partial}(b,w)}
is the complex number of length given by twice the distance of the center of 
the circle circumscribed around (either) face
containing $b,w$ to the mid-point between $b$ and $w$, and direction pointing from
$w$ to $b$.
This is the total Kasteleyn matrix, which we denote by
$$\left(\begin{array}{rr}0& K\\
-K & 0
\end{array}
\right)$$
where the top (resp. bottom) set of rows corresponds to white (resp. black)
vertices. We refer to $K$ as the {\em (discrete) chiral Kasteleyn matrix}.
(In \rref{eanti}, our convention differs from \cite{kenyon}, where the matrix is symmetrical;
however, antisymmetric matrices fit better with Pfaffians.)

In view of the discrete field theory
of Cimasoni-Reshetikhin, the material of Section \ref{fermion} above
suggests a formula of the following form:
Let $\Gamma$ be a graph critically embedded into $\C/\Z\{1,\tau\}$ such
that $[0,1]$ is antiperiodic. Let $b$ be a black vertex and let $f_b(z)$
be the corresponding discrete exponential. Then one should have for the chiral 
Kasteleyn matrix $K$,
\beg{ed*}{Det(K)\sim \cform{\prod}{z:f_{b+1}(z)=-f_b(z)}{} (1\pm f_{b+\tau}(z))
}
where the sign is the spin structure on $[0,\tau]$.

This formula is further motivated by the result of Mercat \cite{mercat} that
exponentials generate additively the vector space of discrete holomorphic
functions; rewritten in the basis of discrete exponentials, the chiral
Kasteleyn matrix becomes diagonal, so its determinant becomes the product
of the diagonal terms, analogously to formula \rref{ed*} in Section \ref{fermion}
above. There are however
two caveats to applying Mercat's theorem: first of all, Mercat is
considering a slightly different setup, not assuming that the graph $\Gamma$
is bipartite. He works with tilings by rhomboids; the graph made by the edges
of the rhomboids is naturally bipartite, with the vertices of the graphs $\Gamma$, $\Gamma^*$
being the vertices of the two different colors. Second, Mercat's result applies
to tilings of simply connected subsets of $\C$, and he does not specify how many
discrete exponentials are needed to get a basis of the space of discrete holomorphic
functions; in order for the formula \rref{ed*} to be strictly correct, we would
have to prove that discrete exponentials $f_{?}(z)$ for such $z$ that
\beg{egper}{f_{b+1}(z)=-f_b(z)} 
precisely form a basis of discrete holomorphic functions
on $\Gamma$ tiling $A_{q}$ twisted by the antiperiodic spin structure on
$[0,1]$ , $q=2\pi i \tau$. One precise statement (for the case of
a rhomboid graph) is as follows:

\vspace{3mm}

\begin{proposition}
\label{trace}
Consider a bipartite graph $\Gamma$ critically embedded into the torus
$E_\tau=\C/\langle 1,\tau\rangle$ ($Im(\tau)>0$) in such a way that there
exist complex numbers $a_1,...,a_{2n}$, $b_1,...,b_{2m}$ of equal absolute
value such that $a_1+...+a_{2n}=1$, $b_1+...+b_{2m}=\tau$, and the corresponding
rhomboid graph on $\Gamma\amalg \Gamma^*$ has edges 
$$(a_1+...+a_{k}+b_1+...+b_\ell, a_1+...+a_{k+1}+b_1+...+b_\ell),$$
$$(a_1+...+a_{k}+b_1+...+b_\ell, a_1+...+a_{k}+b_1+...+b_{\ell+1})$$
(here we consider the subscripts as elements of $\Z/m$, $\Z/n$). Assume further
that the vertex $0$ is black, and assume that the spin structure of $E_\tau$
around $[0,1]$ is antiperiodic. Then the Kasteleyn matrix of $\Gamma$ satisfies
\beg{ekast}{pf(K_\Gamma)=2^m\cform{\prod}{j=1}{n}(1\mp f_\tau(z_j))
}
where $f_b(z)$ are Kenyon exponentials \cite{kenyon}, and $z_j$ ranges
over all complex numbers for which
\beg{eantip}{f_1(z_j)=-1,}
and the sign depends on the spin structure of $[0,\tau]$. 
\end{proposition}

\vspace{3mm}
{\bf Remark:} This is not as easy to generalize as one may think. In particular,
in case of more general Riemann surfaces, a straightforward generalization of this
formula is false. A heuristic explanation may be extracted from Section \ref{fermion}
above: In general, we are unable to identify functions with $1/2$-forms in boundary
behavior, and transformation rules for $1/2$-forms would have to be modelled 
discretely.

\vspace{3mm}
\Proof
Because of continuity, it suffices to prove the formula for generic values
of $a_i$, $b_j$. Consider for $k=0,...,m-1$ the sets of black vertices
of $\Gamma$
$$B_k=\{b_1+...+b_{2k}+a_1+...+a+{2i}\;|\;i=0,...,n-1\}.$$
and the sets of white vertices of $\Gamma$
$$W_k=\{b_1+...+b_{2k+1}+a_{1}+...a_{2i+1}\;|\;i=0,...,n-1\}.$$
The Kenyon exponential function $f_b(z)$ with base point $0$
is on the black vertex
$$b=a_1+...+a_{2\ell}+b_{1}+...+b_{2k}$$
equal to 
$$\cform{\prod}{i=1}{k}\fracd{1+zb_{2i-1}}{1-zb_{2i}}
\cform{\prod}{j=1}{\ell}\fracd{1+za_{2j-1}}{1-za_{2j}}.$$
Recall that we denoted by $z_1,...,z_n$ all numbers satisfying
\rref{eantip}. Consider functions
$$f_{?,k}(z_j):B\r\C$$
where
$$\begin{array}{llll}
f_{b,k}(z_j) & = & f_b(z_j)  &\text{when $b\in B_k$},\\
& = &0 & \text{when $b\notin B_k$.}
\end{array}
$$
We claim that for generic values of $a_i$, $b_j$, $f_{?,k}(z_j)$
form a basis of the space of all functions $B\r \C$. In effect,
it is true for $a_i=\frac{1}{2n}$ by direct computation, and
hence it is true generically (since non-degeneracy occurs on a 
Zariski-open set).

Now by direct computation,
\beg{edebar+}{
\begin{array}{l}
\overline{\partial}f_{?,k}(z_j)=g_{?,k}(z_j)-g_{?,k-1}(z_j)\;
\text{for $k=1,...,m-1$}\\
\overline{\partial}f_{?,0}(z_j)=g_{?,0}(z_j)\mp f_\tau(z_j)g_{?,n}(z_j)
\end{array}
}
where 
$$g_{?,k}(z_j):W\r\C$$
are functions defined as follows:
\beg{edefg}{
\begin{array}{l}
g_{b_1+...+b_{2\ell+1}+a_1+...+a_{2i+1},k}(z_j):=\\
f_{b_1+...+b_{2\ell}+a_1+...+a_{2i+2},k}(z_j)-
f_{b_1+...+b_{2\ell}+a_1+...+a_{2i},k}(z_j).
\end{array}
}
Note that for each $j$, the matrix \rref{edebar+} in the basis
$(f_?,k(z_j),g_{?,k}(z_j))$ is
$$
\left(
\begin{array}{crrrrr}
1&-1&&&&\\
&1&-1&&&\\
&&1&-1&&\\
...&...&...&...&...&...\\
\mp f_\tau(z_j)&&&&&1
\end{array}
\right),
$$
which has determinant 
$$1\mp f_\tau(z_j).$$
But now the function $h_{?,k}:W\r \C$ which is defined by
$$h_{b_1+...+b_{2\ell+1}+a_1+...+a_{2i+1},k}(z_j):=
f_{b_1+...+b_{2\ell}+a_1+...+a_{2i+2},k}(z_j)$$
plays a symmetric role to $f_{?,k}$ which respect to exchanging
$\overline{\partial}$ for $\partial$, and the base change
between $g$ and $h$ is conjugate to the matrix
$$
\left(
\begin{array}{rrrrrr}
1&-1&&&&\\
&1&-1&&&\\
&&1&-1&&\\
...&...&...&...&...&...\\
1&&&&&1
\end{array}
\right),
$$
which is $2$.
\qed

\vspace{3mm}
To try to obtain a limit over
tilings with decreasing length of rhomboid edge, it further seems natural
to use a recipe of Mercat \cite{mercat}
how to take the limit in order to make discrete exponentials
converge to ordinary continuous exponentials: subdivide each rhomboid into 4 equal
rhomboids with parallel edges, and repeat this procedure. Even in simple examples, however, we see that
this procedure when applied to \rref{ed*}, will {\em not} produce a convergent limit
in the naive sense: the limit has to be ``regularized'' (which is not surprising,
given the usual physical context of lattice regularization). For example, 
if the segment $[0,1]$ is tiled into equal segments by the edges of $\Gamma$,
the equation \rref{egper} becomes
$$\left(\frac{1-z/N}{1+z/N}\right)^{2N}=-1,$$
which implies that $z$ is on the imaginary line. Then \rref{ed*} becomes a product
of quantities of the form
\beg{egper1}{1\pm e^{-z\tau},}
so given $Im(\tau)>0$, for $Im(z)>0$, the second summand
\rref{egper1} will have absolute value $>1$.
The number of such terms will increase
with $N$ , as solutions of \rref{egper} in increasing
bounded regions will
approach $(2k+1)\pi i$, $k\in \Z$. Therefore, the expression \rref{ed*} in this
case cannot converge as $N\r \infty$. The ``limit'' has to be taken in the sense
of 
$$\cform{\prod}{n\in \Z}{}(1-q^{n+\frac{1}{2}})=
\pm \cform{\prod}{n\in \N}{}(1-q^{n+\frac{1}{2}})^2 q^{-\frac{1}{2}-\frac{3}{2}-\frac{5}{2}-...}\sim$$
$$\sim\cform{\prod}{n\in\N}{}(1-q^{n+\frac{1}{2}})^2 q^{-\frac{1}{24}},$$
(the right hand side being the desired ratio of a theta function by $\eta$), 
based on the ``calculation"
$$-\frac{1}{2}-\frac{3}{2}-\frac{5}{2}-...\sim$$
$$\sim(-\frac{1}{2}-1-\frac{3}{2}-2-\frac{5}{2}-...)
+(1+2+3+...)\sim -\frac{1}{2}\zeta(-1)+\zeta(-1)=-\frac{1}{24}.$$
In this special case, the regularization can be achieved by multiplying \rref{ed*}
by
\beg{egper2}{\cform{\prod}{Im(z)>0}{}\prod \frac{1}{f_{b+\tau}(z)}}
where the product is over solutions of \rref{egper} with $Im(z)>0$. This regularization
indeed produces the desired ratio of $\theta$ and $\eta$ in the limit, multiplied
by the factor $q^{1/24}$. 

Therefore, in some sense, for a general graph, the question entails
finding a regularization procedure which would generalize \rref{egper2}.
Again, one precise statement can be obtained as follows. Define
$$:pd (K_\Gamma): \;= \cform{\prod}{j=1}{n}(1-f_\tau(i|z_j|)).$$

\vspace{3mm}

\begin{proposition}
\label{plimit}
Suppose graphs $\Gamma_k$ are as in Proposition \ref{trace}, and that
the rhomboid graph of
$\Gamma_{k+1}$ is obtained from the rhomboid graph of $\Gamma_k$
by subdividing each rhomboid into 4 congruent rhomboids with parallel
edges. Then
$$\cform{\lim}{k\r\infty}{}:pf (K_{\Gamma_k}): =\cform{\prod}{j=0}{\infty}
(1\mp q^{j+\frac{1}{2}})^2.$$
\end{proposition}

\Proof
We must study solutions $z_j$ of the equation
$$\cform{\prod}{i=1}{n} \fracd{\left(1+\fracd{a_{2i-1}z}{N}\right)^N}{
\left(1-\fracd{a_{2i}z}{N}\right)^N}=-1.$$
We rewrite this as
\beg{elim1}{\cform{\sum}{i=1}{n}\left(
\ln\left(1+\fracd{a_{2i-1}z}{N}
\right)
-
\ln\left(1-\fracd{a_{2i}z}{N}
\right)
\right)
=\fracd{2k+1}{N}\pi i,\; k\in\Z.
}
Put
$$f(t):=
\cform{\sum}{i=1}{n}(\ln(1+a_{2i-1}t)-\ln(1-a_{2i}t)).
$$
Then
$$f(0)=0$$
and $f$ is analytic in a neighborhood of $0$ with non-zero derivative at
$0$. Hence, the same holds for the inverse $g=f^{-1}$. From this,
(substituting $t=\frac{z}{N}$), we have the following

\vspace{3mm}
\begin{lemma}
\label{l1}
For every $\delta>0$ there exists a $c>0$ such that for every $N$, when
\beg{elim+}{|z|<cN,
}
then $\frac{z}{N}$ is in the domain of $g$, and all the solutions to
\rref{elim1} satisfying \rref{elim+} are of the form
\beg{elim++}{z_k=(2k+1)\pi i+\epsilon_k, \; |\epsilon_k|<\frac{\delta k}{N}.
}
\end{lemma}
\qed

\vspace{3mm}

Now consider an inequality of the form
\beg{elim2}{C^{-1/M}<\left|
\fracd{1-q^{|z|/2\pi}}{1-q^{|k+1/2|}}
\right|<C^{1/M}
}
where $C>1$ is some constant. \rref{elim2} holds when
$$\left|
\ln|1-q^{|z|/2\pi}|-\ln|1-q^{|k+1/2|}|
\right|<\frac{\ln C}{M},$$
or
$$|q^{|z|/2\pi}-q^{|k+1/2|}|<\frac{\ln C}{2M}.$$
This follows whenever
\beg{elim3}{\left|
\frac{|z|}{2\pi}-(|k+\frac{1}{2}|)
\right|<\frac{k}{M}
}
for some constant $K>0$ dependent on $C$, since $q^t$ has a bounded derivative
in $t>0$ for $0<|q|<1$.

Now set $M:=\sqrt{N}$. By the Lemma, there exist constants $A,B>0$, 
(we may assume $A<1$) such that \rref{elim3} (and hence \rref{elim2})
holds for some solution $z_j$ whenever
\beg{elim4}{|k|<A\sqrt{N},
}
and
\beg{elim5}{|z_j|>B\sqrt{N}\;\text{for all $j\in S$}
}
where $S$ is the set of all other indices $j$. Now note that
\beg{elim6}{\cform{\prod}{j\in S}{}(1-q^{|z_j|/2\pi})=
O(\cform{\sum}{j\in S}{}(|q^{|z_j|/2\pi}|)
=O(N|q|^{B\sqrt{N}/2\pi}).
}
The $\limsup$ of the right hand side of \rref{elim6} with $N\r\infty$ is $0$.
Thus, the solutions \rref{elim5} may be neglected. Also, the quantity
$$\cform{\prod}{k>A\sqrt{N}}{}(1-q^{|k+1/2|})^2$$
approaches $0$ as $N\r\infty$, and by \rref{elim2},
$$\frac{1}{C}<\cform{\prod}{|k|<A\sqrt{N}}{}\left|
\fracd{1-q^{|z|/2\pi}}{1-q^{(|k+1/2|)}}
\right|<C,$$
as $M=\sqrt{N}$ and $A<1$. Since $C>1$ was an arbitrary constant,
we are done.
\qed


\vspace{3mm}
\section{Fermions on a Riemann surface II: A brief review of conformal field theory,
with speculations about graphs}

\label{fcft}

In physics, it is ``known'' that 
the critical Ising model on a lattice converges
to the $1$-dimensional free fermion theory \cite{onsager}, cf. \cite{difr},
while the critical dimer model converges to the free fermion
of space-time dimension 2 (\cite{dijk}). 
To obtain rigorous mathematical
theorems in these directions, 
we need to take into account the detailed structure
of the fermionic CFT's, and model it in graph theory. Such
mathematical treatments of conformal field theory now exist,
(cf. \cite{scft,spin} from the Segal point of view, \cite{feingold}
from the vertex algebra point of view),
but the structure involved
is somewhat complicated. Because of this,
we take time in this section to survery these structures in more
detail. Some comments on their potential impact on the graph theory side
of the story will be made at this section's conclusion. For simplicity, 
throughout this section, we will consider {\em complex} CFT structures,
as real structures complicate things further.

\vspace{3mm}
A big part of the story of fermion CFT's is {\em bosonisation}.
For a bosonic CFT, we are supposed to have a real-analytic line bundle $L$
on the moduli space of Riemann surfaces $\Sigma$ with analytically parametrized
boundary, a complex Hilbert space $H$ and, writing $H_{\partial\Sigma}$ for
$$H^{*\hat{\otimes}m}\hat{\otimes} H^{\hat{\otimes}n}$$
when $\Sigma$ has $m$ inbound and $n$ outbound boundary
components (we take $\C$ when $\Sigma$ is closed), a map
\beg{ecft*}{L_\Sigma\r H_{\partial \Sigma}
}
satisfying gluing axioms analogous to Lemma \ref{lg}
(cf. \cite{scft,hk,fhk}). 

\vspace{3mm}
A genus $g$ Riemann surface $\Sigma_g$ can be cut along $g$ non-separating disjoint
real-analytic simple curves $a_1,...,a_g$ into a genus $0$ surface $A$ with $2g$
boundary components (let, for future reference, a ``dual cutting'' by real-analytic
simple curves $b_1,...,b_g$.). Gluing a standard disk on each of these boundary components,
we obtain a sphere. Since the $2$-sphere possesses a unique
conformal structure, this specifies an element
$$\iota_{a_1,...,a_g}\in L_{\Sigma_g},$$
hence, by \rref{ecft*}, a complex number. This assignment
\beg{ehigherp}{A\mapsto Z(a_1,...,a_g)\in \C}
is called the {\em genus $g$ partition function}. For $g=1$,
one can take the standard annulus $A_q\subset \C$ with boundary components
$S^1$ and $qS^1$ parametrized by the functions $z$, $qz$, respectively, which
gives us the {\em partition function} $Z(q)$. This function is a series 
in
\beg{eweight}{q^{a}\overline{q}^b} 
($0<||q||<1$) with positive integral
coefficients where $a-b\in\Z$, and the coefficient of \rref{eweight} counts
the number of {\em states of $H$ of weight $(a,b)$}. In other words, it
is the trace of the ``grading operator'' which, at weight $(a,b)$,
is given by multiplication by $q^{a}\overline{q}^b$.

\vspace{3mm}
There is a special class of CFT's called {\em rational conformal field theory},
abbr. RCFT. (One rigorous mathematical approach was reached by Huang and Lepowsky
\cite{hl}, although it doesn't use the language of Hilbert spaces, but
the structure of vertex algebra. An approach
via Hilbert spaces was developed in \cite{hk,fhk}. To date, the two approaches
have not been fully unified, and each has certain desirable features.) In this case,
we have a finite set $\Lambda$ of {\em sectors} with a distinguished element
$0\in\Lambda$ and an involution $(?)^*:\Lambda\r \Lambda$, $0^*=0$, and
decompositions
$$H=\bigotimes_{\lambda\in\Lambda} H_{\lambda}\hat{\otimes} \overline{H}_{\lambda^*},$$
\beg{ecftfactor}{Z(q)=\sum_{\lambda\in\Lambda}Z_{\lambda}(q)Z_{\lambda^*}(\overline{q}).}
(Sometimes, elements of $\Lambda$ are referred to as {\em labels}, an indexing
set for sectors. We choose here to identify each label with its corresponding
sector.)
The precise explanation of \rref{ecftfactor} is that we have certain
line bundles $M_\lambda$ on the moduli space of all 
annuli, and a positive-definite
Hermitian pairing
\beg{epair}{\bigoplus_{\lambda\in\Lambda} M_{\lambda}\otimes \overline{M}_{\lambda^*}
\r L,
}
and we have
$$Z_{\lambda}(q):M_\lambda A_q\r\C$$
so \rref{ecftfactor} makes sense. 

This is actually the tip of an iceberg called {\em chiral conformal field theory}:
One is required to have a finite-dimensional complex vector space $M_\lambda\Sigma$
for each Riemann surface $\Sigma_{(\lambda_i)}$ with parametrized boundary components
labeled by the $\lambda_i$'s, which are supposed
to satisfy both \rref{epair} and gluing axioms \cite{hk, fhk}. The most important
part of the gluing structure is an isomorphism
\beg{emodglue}{M{\check{\Sigma}_{(\lambda_i)}}\cong \bigoplus_{\lambda\in\Lambda} 
M\Sigma_{(\lambda_i,\lambda,
\lambda)}}
where $\check{\Sigma}$ is obtained from $\Sigma$ by gluing an inbound and
an outbound boundary component: the sum on the right hand side is over all possible
labellings of the two additional boundary components by the same label.
Additionally, there
are to be linear maps 
\beg{echiralop}{M\Sigma_{(\lambda_i)}\r H_{\partial \Sigma_{(\lambda_i)}}}
satisfying appropriate gluing axioms (\cite{hk,fhk}). It is worth mentioning
that for each $\lambda\in\Lambda$, the line bundle $M_\lambda A$ specifies
a holomorphic $\C^\times$-central extension of the semigroup of annuli with
analytically parametrized boundary components of opposite orientations; these
central extensions do not depend on $\lambda$,
and are characterized by a single number called the {\em central
charge} and denoted by $c$ - see \cite{scft}. (Note: there is
an infinite-dimensional space of annuli with analytically parametrized
boundary components modulo the relation of conformal equivalence preserving
the boundary parametrizations; if we restricted only to the standard
annuli $A_q$, we would get the semigroup $\C_{0<q<1}$ which does not
have non-trivial central extensions.)

\vspace{3mm}

Models of the fermion CFT which satisfy the above axioms are referred to as
{\em bosonised fermions}. From the point of view of string theory, 
a fermion of spacetime dimension $d$ has central
charge $d/2$. The theory differs depending
on whether $d$ is odd or even. For $d=1$, the theory
can be constructed, for example, as a coset model of the level $2$ WZW-model
for $SU(2)$ (\cite{difr}, 18.5.1, $k=2$), or as a non-supersymmetric minimal model 
\cite{difr}, 7.4.2.

\vspace{3mm}
Generally, for $d$ odd, the $d$-dimensional bosonised fermion theory
has three sectors 
$$\Lambda=\{NS^+,NS^-,R\},\;\lambda^*=\lambda,\; 0=NS^+.$$
(R stands for Ramond and NS stands for Neveu-Schwarz.)
Its partition function is 
determined by letting $Z_{NS^+}(q)$, $Z_{NS^-}(q)$ be the integral (resp. 
integral $+1/2$)-dimensional terms of
$$\prod_{n\geq 0} (1+q^{n+\frac{1}{2}})^{d},$$
and 
$$Z_R(q)=2^{\lfloor d/2\rfloor}\prod_{n>0} (1+q^n).$$
The power $2$ in the formula for $Z_R(q)$ is significant, it
expresses the fact that the bottom weight space of $H_R$ is naturally
the spinor module of the complex Clifford algebra on a space of dimension $d$.
The dimension of the modular functor on a Riemann surface of genus $g$ is
\beg{eodddim}{\frac{1}{2}(4^g+2^g),
}
which is equal to the number of spin structures of Arf invariant $0$.

\vspace{3mm}
The bosonised fermion of even spacetime dimension $d$ 
has sector set
$$\Lambda_2=
\{NS^+,
NS^-, R_+, R_-\}.$$
We have
$$R_{+}^{*}=R_-.$$
Here the labels we chose express the operations on the sector Hilbert spaces
of the $1$-dimensional fermion. In particular, $H_R\hat{\otimes} H_R$ 
decomposes into two state spaces, according to the decomposition of the
complex spinor on $\R^2$.
Its partition function is 
determined by letting $Z_{NS^+}(q)$, $Z_{NS^-}(q)$ be the integral (resp. 
integral $+1/2$)-dimensional terms of
$$\prod_{n\geq 0} (1+q^{n+\frac{1}{2}})^{d},$$
and 
$$Z_{R_\pm}(q)=2^{ d/2}\prod_{n>0} (1+q^n).$$
The reason why there are now two Ramond sectors is that there are
two non-isomorphic complex spinors in even dimension.
The dimension of the modular functor on a Riemann surface of genus $g$ is
\beg{eevendim}{4^g,
}
which is equal to the number of all spin structures.

\vspace{3mm}
For $d=2$, the bosonised fermion theory turns out to be {\em isomorphic} to the 
conformal field theory associated with the lattice $2\Z\subset \R$ (\cite{fk, hk}).
This phenomenon is known as the {\em boson-fermion correspondence}.
Expressed this way, the labels become
$$\{0,1/2,1,3/2\},\; (1/2)^*=(3/2).$$
The chiral partition functions are
\beg{elattice}{Z_{\lambda}(q)=\sum_{k\in\Z} q^{\frac{1}{2}(k+\lambda)^2}
\prod_{n>0}(1-q^n)^{-1}.}

\vspace{3mm}
We see that regardless of the parity of dimension, the partition function of
the fermion of dimension $d$ is
\beg{efpart}{\begin{array}{l}
\displaystyle
Z(q)=\frac{1}{2}\left(\prod_{n\geq 0}(1+q^{n+1/2})(1+\overline{q}^{n+1/2})+
\prod_{n\geq 0}(1-q^{n-1/2})(1-q^{n-1/2})+\right.\\[4ex]
\displaystyle
\left. 2^{d+1}(q\overline{q})^{d/16}
\prod_{n\geq 1}(1+q^{n})(1+\overline{q}^{n})\right).
\end{array}}
Note that for $d=2$, this is indeed equal to 
\beg{ed2part}{\begin{array}{l}
\left(\sum_{k\in\Z}((q\overline q)^{2k^2}+(q\overline q)^{\frac{1}{2}(2k+1)^2}
+ 2(q\overline q)^{\frac{1}{2}(2k+\frac{1}{2})^2})\right)\cdot\\[4ex]
\displaystyle \cdot \prod_{m,n>0} (1-q^n)^{-1}(1-\overline{q}^n)^{-1},
\end{array}
}
which is what we obtain when we combine the chiral partition functions
\rref{elattice}.
Note that in \rref{ed2part} and \rref{efpart}, in fact the sums of the 
first two summands are equal by Jacobi's triple product 
identity, as are the last summands.

\vspace{3mm}
How is it possible that the dimensions of the modular functor \rref{eodddim} and 
\rref{eevendim} in the even and odd-dimensional case are different, even though
the partition function \rref{efpart} has a uniform expression? The answer is
that in the even-dimensional case, when the spin structure along some of the curves
$a_i$ is anti-periodic, the two Ramond sectors arising will have the same
partition function, and hence their sum can be grouped in the power of $2$
coefficient.

\vspace{3mm}
The {\em genuinely fermionic} approach to fermion RCFT uses another,
more complicated, axiomatization.
This time, the chiral parts are {\em chiral super conformal field theories}. 
(This is a different concept from a supersymmetric conformal
field theory, where both a bosonic and a fermionic part
are present.) Chiral super conformal field
theories are axiomatized similarly as chiral conformal
field theories, but the Riemann surfaces $\Sigma$ involved come
with spin structure. We then have a decomposition
$$\Lambda=\Lambda_{NS}\amalg \Lambda_{R},$$
expressing which labels apply to antiperiodic and which to periodic boundary components. 
Additionally, we deal with {\em super-modular functors}, which means that each
$M\Sigma_{(\lambda_i)}$ is a $\Z/2$-graded vector space, as is each
sector state space
$$H_{\lambda}.$$
The maps \rref{echiralop} are required to be {\em graded}, and
the usual sign convention 
$$(-1)^{ij}$$
is applied with interchanging elements of degree $i$, $j$. We now have 
chiral partition functions
$$Z_{\lambda}^{NS}(q),\; Z_{\lambda}^{R}(q)$$
where the superscript expresses the spin structure in the radial direction
of the annulus. When factoring a CFT to fermionic chiral 
factors, the analog of formula \rref{ecftfactor} at both chiralities
becomes
\beg{ecftfactori}{Z(q)=\frac{1}{2}\left(\sum_{\lambda\in\Lambda} Z^{NS}_{\lambda}(q)
Z^{NS}_{\lambda^*}(\overline{q})
+\sum_{\lambda\in\Lambda} Z^{R}_{\lambda}(q)
Z^{R}_{\lambda^*}(\overline{q})\right).}
For a Riemann surface of
arbitrary genus $g$ where the summation is over spin structures on the 
Riemann surface $\Sigma$,
and over a product of the sets of periodic (resp. antiperiodic) labels matching
to the spin structure of each of the curves $a_i$, $b_i$, and the
coefficient before the sum is $1/2^g$. (This can be thought of
as an averaging over spin structures, or,
in analogy with string theory, a ``GSO projection''.) When
there is exactly one periodic and one anti-periodic label, 
one can show that the super-modular functor is always one-dimensional
for each spin structure on a Riemann surface, and is {\em invertible}
in the sense that its tensor product with another super-modular functor
is the trivial (unit) super-modular functor.
In this case, the summation is {\em precisely} over the set of spin structures
on $\Sigma$. This, in fact, is the CFT explanation of
the decomposition of the discrete partition functions of the Ising and dimer
models into summands corresponding to spin structures, with signs given
by the Arf invariant (see e.g. formula \ref{egen} above); those are the summands 
whose gluing we study
in Sections \ref{s.dg} and \ref{sec.str}.
There is one substantial difference in the ``CFT limit'', however:
one can show that {\em the summands corresponding
to spin structures of Arf invariant $1$ are always equal to $0$}.
From a physics point of view, this can be justified by noting that
the contribution of such terms would come with a negative sign, which
is physically impossible. A mathematical argument can also be obtained
along the lines of \cite{spin}.

\vspace{3mm}
It is also possible to talk about a fermionic conformal field theory
with both chiralities, which consists of two super-Hilbert spaces $H_R$ and
$H_{NS}$ ($R$ corresponds to the periodic and $NS$ to the anti-periodic spin structure
on $S^1$), a real-analytic line bundle $L$ on the moduli space of Riemann surfaces
with spin structure and with
(analytically) parametrized boundary components. The conformal field theory
then specifies, again, a map
$$L_\Sigma\r H_{\partial \Sigma}$$
where $H_{\partial\Sigma}$ is, again, a Hilbert tensor product of copies of
$H_R$, $H_{NS}$ and their duals, depending on orientation and spin structure
of each boundary components. The gluing axioms are analogous to the bosonic case
(Note: no super-structure is present on $L$, or in other words, it is
assumed to be ``even''.)

Given a fermionic chiral RCFT, the spaces $H_{R}$ and $H_{NS}$ are
given by
$$\bigoplus_\lambda H_{\lambda}\hat{\otimes} \overline{H}_{\lambda^{*}}$$
where the sum is taken over R resp. NS labels. For a closed Riemann surface
$\Sigma$, one then has one partition function for each spin structure, although
the partition functions corresponding to spin structures of Arf invariant $1$ are
$0$. The sum of these fermionic partition functions, multiplied by the factor $1/2^g$,
is the bosonic partition function \rref{ecftfactori}. 

\vspace{3mm}
There is a {\em $d$-dimensional} chiral
free fermion with invertible super-modular functor
when {\em $d$ is even}. The chiral fermionic partition functions are
\beg{eferm1}{{}_cZ_{NS}^{R}(q)=\prod_{n\geq 0}(1+q^{n+\frac{1}{2}})^d,}
\beg{eferm2}{{}_cZ_{NS}^{NS}(q)=\prod_{n\geq 0}(1-q^{n+\frac{1}{2}})^d,}
\beg{eferm3}{{}_cZ_{R}^{NS}(q)=2^{d/2}\prod_{n> 0}(1+q^{n})^d,}
\beg{eferm4}{{}_cZ_{R}^{R}(q)=0}
(We use the symbol ${}_cZ$ to distinguish the chiral partition functions from
the fermionic partition functions at both chiralities, which are written
below in \rref{efermp1}, \rref{efermp2}, \rref{efermp3}, \rref{efermp4}, since
there is only one label for each spin structure.)
The conformal field theory assembled from this super-RCFT is
isomorphic to the conformal field theory
of the bosonised $d$-dimensional fermion described above.

\vspace{3mm}
For $d$ odd, {\em there is no $d$-dimensional chiral fermion
super-conformal field theory} as defined above.
It was proved in \cite{spin} that there is no invertible super-modular functor of central
charge $1/2$. However, P.Deligne found an even more involved axiomatization using
the {\em super Brauer group} and central simple algebras which does allow 
a $1$-dimensional chiral fermion model. Roughly, the super-Brauer group of 
$\C$ is $\Z/2$, and then non-zero element is represented by the super-division
algebra $C=\C[a]/(a^2-1)$ where $a$ is odd. Then the value of the modular functor
(and the Hilbert spaces) are taken in the category of modules over super-division
algebras, and identification is allowed under super-Morita equivalence. 
When $d$ is odd, $1$-dimensional $C$-modules occur as values of the modular
functor when Ramond boundary components arise. 

\vspace{3mm}
When we put both chiralities together, there is, again, a uniform formula
for the fermionic (spin-structure dependent)
partition functions of the $d$-dimensional fermion, regardless of the parity
of $d$:
\beg{efermp1}{Z_{NS}^{R}(q)=\prod_{n\geq 0}(1+q^{n+\frac{1}{2}})^d(1+\overline{q
}^{n+\frac{1}{2}})^d,}
\beg{efermp2}{Z_{NS}^{NS}(q)=\prod_{n\geq 0}(1-q^{n+\frac{1}{2}})^d
(1-\overline{q}^{n+\frac{1}{2}})^d,}
\beg{efermp3}{Z_{R}^{NS}(q)=2^{d}\prod_{n> 0}(1+q^{n})^d(1+\overline{q}^{n})^d,}
\beg{efermp4}{Z_{R}^{R}(q)=0}
In fact, as a fermionic conformal field theory (i.e. on spin structure at a time),
the $d$-dimensional fermion is a tensor product of $d$ copies of the $1$-dimensional
fermion (and its rigorous mathematical construction in fact proceeds in this fashion).

For $d$ odd, the Hermitian form 
$$C\otimes \overline{C}\r \C$$
has
$$\langle 1,1\rangle=\langle a,a\rangle=1,$$
which causes, on the R label part of the modular functor, which is given by
Deligne's prescription in the
form of a (1-dimensional) graded $C$-module, the additional factor $2$, explaining
the odd power of $2$ present in \rref{efermp3}.

\vspace{3mm}

\vspace{5mm}
{\bf Comments and questions about graphs:}
Physics \cite{dijk} predicts that there 
should exist a limit formula from the critical dimer model to the $2$-dimensional
free fermion conformal field theory, i.e. the theory which has a bosonic
description with modular functor with set of labels $\Lambda_2$, as well as a fully
fermionic description with an invertible super-modular functor.

On the other hand, physics also predicts \cite{onsager} that there should exist
a limit formula from the critical Ising model to the bosonic $1$-dimensional 
fermion theory with the set of labels $\Lambda$. We saw that both conformal
field theories are closely related yet substantially different. 

The factorization formulas \rref{ecftfactor} and \rref{ecftfactori},
and their higher genus analogues,
suggest compelling analogies with decompositions of partition functions
of the discrete models into summands. 
This is particularly striking
in view of the expression of the dimer and Ising model partition functions
as a sum of graphs drawn in a suitable way on genus $g$ Riemann
surfaces as sums of $4^g$ Pfaffians. In fact, an important point
of this paper is that in Sections \ref{s.dg} and
\ref{sec.str}, we have modeled discrete gluing analogues for
the {\em fermionic} form of the conformal field theory, i.e. \rref{efermp1},
\rref{efermp2}, \rref{efermp3}, \rref{efermp4}.

It is not clear whether (or in what situations)
the estimates for number of summands given by the limit theory are
precise: In the conjectured limit, in the fermionic description, \rref{efermp4},
and more generally, 
the partition functions corresponding to spin structures of Arf
invariant $1$, always vanish (see also \cite{dijk}).
This may appear to contradict the result of Norin \cite{norin} showing that 
for $g=1$, a $3$-Pfaffian graph is always $1$-Pfaffian,
while the limit formula predicts $3$ non-vanishing factors. However, Norin's
result \cite{norin} is stated in a form not involving weights of edges.
When translated to our formulation with weights, cancellations
may occur, and some terms may vanish. This
occurs specifically in the case of critical embeddings, even without taking the
limit.

Note that in the bosonic description, the interpretation of the limit
is somewhat different. In this case, the summands are not Pfaffians, but
suitable averages of Pfaffians using the GSO projection. In this case,
in the even-dimensional case, none of the summands actually vanish,
but subsets of different summands have identical partition functions:
in the discrete case, then, the question of vanishing is replaced by other linear
dependencies.

In the odd-dimensional case (which is relevant to the Ising model),
the situation is more interesting: In the bosonic description, 
one doesn't combine sectors with the same partition function,
but instead the sectors themselves combine to one sector. Is there
a definition of corresponding ``GSO-projected" summands purely in
terms of graphs?

Alternately, in the fermionic description, is there a graph-theoretic
characterization of Pfaffians corresponding to the modules over the
Deligne super-algebra $C$? While an expression of the Ising partition
function as a sum of $4^g$ Pfaffians is obtained in \cite{LM}, this 
description uses a modifying construction on the underlying graph.
Is there a formula in terms purely of Pfaffian-like summands on the
original graph, suitably generalized using
$C$-modules discretely mimicking Deligne's theory?


\vspace{10mm}
\begin{thebibliography}{99}


\bibitem{agmv} L. Alvarez-Gaume, G. Moore, C. Vafa: 
Theta functions, modular invariance, and strings, Comm. Math. Phys. Volume 106, Number 1 (1986), 1-40.

\bibitem{bd} A. A. Beilinson, V. Drinfeld: Chiral algebras, ISBN-10: 0-8218-3528-9, AMS Colloquium Publications, vol. 51.

\bibitem{bor} R. E. Borcherds: Monstrous Moonshine and Monstrous Lie Superalgebras, Invent. Math. 109, 405�444, 1992.


\bibitem{cr1} Cimasoni, D., Reshetikhin, N.: Dimers on surface graphs and spin structures, 
Communications in Mathematical Physics, Vol. 275, 2007, p. 187-208.

\bibitem{cr} Cimasoni, D., Reshetikhin, N.: Dimers on surface graphs and spin structures II., 
Communications in Mathematical Physics, Vol. 281, No. 2, 2009, p. 445-468.

\bibitem{ckp} H.Cohn, R. Kenyon, J.Propp: A variational principle
for domino tilings, {\em J. Amer. Math. Soc.} 14 (2001) 297-346

\bibitem{difr} P. Di Francesco, P. Mathieu, and D. S\'{e}n\'{e}chal: 
Conformal Field Theory, Springer-Verlag, New York, 1997. ISBN 0-387-94785-X.

\bibitem{dijk} R. Dijkgraaf, D. Orlando, S. Reffert: 
Dimer models, free fermions and super quantum mechanics, 
Adv. Theor. Math. Phys. Volume 13, Number 5 (2009), 1255-1315.

\bibitem{feingold} A.Feingold, F.X.Ries, M.D.Weiner: Spinor construction of
the $c=1/2$ minimal model, in: {\em Moonshine, the Monster and Related Topics},
Contemporary Math. 193 (C.Dong and G. Mason, eds.), AMS, Providence, R.I., 1995, 45-92

\bibitem{ferdinand} J.Ferdinand: Statistical mechanics of dimers on a 
quadratic lattice, {\em J.Math. Phys.} 8 (1967), 2332-2339

\bibitem{fhk} T.Fiore, P.Hu and I.Kriz: 
Laplaza sets, or How to Select Coherence Diagrams for Pseudo Algebras, 
Advances in Mathematics 218 (2008) 1705-1722.



\bibitem{fk} T.Fiore, I.Kriz: What is the Jacobian of a Riemann surface with boundary?,
Deformation spaces,  53�74,
Aspects Math., E40, Vieweg + Teubner, Wiesbaden, 2010

\bibitem{f} M. E. Fischer:
On the dimer solution of planar Ising problems,
Journal of Mathematical Physics 7 10 (1966).

\bibitem{flm} I. Frenkel, J. Lepowski, and A. Meurman: Vertex operator algebras and the monster, 
Academic Press, Boston, 1988.

\bibitem{gl}
A. Galluccio, M. Loebl: A Theory of Pfaffian Orientations I., Electronic Journal of Combinatorics 6,1 1999.


\bibitem{has}
{\sc A. Hassell}, Analytic surgery and analytic torsion, Comm. in Anal. and Geom. 6, 255-289 (1998).

\bibitem{hkp} R.Hortsch, I.Kriz, A.Pultr: A universal approach to vertex algebras, 
J. Algebra. 324 (2010), 1731-1753.


\bibitem{hk} P. Hu, I. Kriz: Conformal field theory and elliptic cohomology, Advances in Mathematics,
Volume 189, Issue 2, 20 December 2004, Pages 325-412. 

\bibitem{hl} Y.-Z. Huang, J. Lepowski: A theory of tensor products for module categories
for a vertex operator algebra. I, II., Selecta Math. (N.S.) 1 (1995),
no. 4, 699�756, 757�786.

\bibitem{kast} P. W. Kasteleyn, Dimer statistics and phase transitions. J. Math. Phys. 4 (1963), 287�293.

\bibitem{kenyon} R. Kenyon: The Laplacian and Dirac operators on critical planar graphs, Invent. Math.,
150 (2002), pp. 409�439.

\bibitem{spin} I. Kriz: On Spin and modularity in conformal field theory, Ann. Sci. de ENS 36 (4) (2003), 1, 57�112.


\bibitem{kx} I.Kriz, Y.Xiu: Tree field algebras: An algebraic axiomatization of
intertwining vertex operators, arXiv:1102.2007

\bibitem{LM}
M. Loebl, G. Masbaum, On the optimality of the Arf invariant formula for the graph polynomials,
 Advances in Mathematics
Volume 226, Issue 1, 15 January 2011, Pages 332-349.

\bibitem{ls}
{\sc M. Loebl, P. Somberg}, Discrete Dirac Operators, Critical Embeddings and Ihara-Selberg Functions,
http://xxx.lanl.gov/abs/0912.3200.


\bibitem{mercat} C. Mercat: Exponentials form a basis of discrete holomorphic
functions, arXiv: math-ph/0210016.


\bibitem{norin} S. Norin: Drawing 4-Pfaffian graphs on the torus, Combinatorica 29 (2009), 109-119.

\bibitem{norin6} A.A.A.Miranda, C.L.Lucchesi: Matching signatures and Pfaffian graphs, 
Discrete Mathematics, 2011, 289-294.


\bibitem{onsager} L. Onsager: Crystal statistics. I. A two-dimensional model with an order-disorder transition, 
1944, Phys. Rev. (2) 65 (3�4), 117�149.

\bibitem{ps} A. N. Pressley and G. B. Segal: Loop groups, Oxford mathematical monographs, 1988, 318 pages.

\bibitem{scft} G. Segal: The definition of conformal field theory, in: Differential geometrical methods in theoretical physics (Como, 1987), NATO Adv. Sci. Inst. Ser. C Math. Phys. Sci., 250, Kluwer Acad. Publ., Dordrecht, 1988, 165-171.

\bibitem{tes}
G. Tesler: Matching in graphs on non-orientable surfaces, J. Comb. Theory B 78
(2000), 198�231.

\bibitem{vdW}
B. L. van der Warden: Die lange Reichweite der regelm�assigen Atomanordnung in Mischkristallen, Z. Physik 118:473,
(1941).

\bibitem{w} E. Witten: Free fermions on an algebraic curve, 
The Mathematical Heritage of Hermann Weyl (Durham, NC, 1987), 
Proc. Sympos. Pure Math., vol. 48, Amer. Math. Soc., Providence, 1988, pp. 329�344.

\end{thebibliography}
\end{document}